\documentclass[12pt]{amsart}
\usepackage{amssymb,times}
\newtheorem{thm}{Theorem}
\newtheorem{prop}{Proposition}
\newtheorem{lem}{Lemma}
\newtheorem{cor}{Corollary}

\newtheorem{sa}[thm]{Theorem}

\newtheorem{hs}[lem]{Lemma}

\theoremstyle{remark}
\newtheorem{rem}{Remark}
\newtheorem{ex}{Example}

\theoremstyle{definition}

\newcommand{\C}{\mathbb{ C}}

\newcommand{\map}{\rightarrow}

\newcommand{\weg}{\smallsetminus}
\newcommand{\eps}{\varepsilon}
\renewcommand{\ge}{\geqslant}
\renewcommand{\geq}{\geqslant}
\renewcommand{\le}{\leqslant}
\renewcommand{\leq}{\leqslant}

\newcommand{\PP}{\mathbb P}

\title[Einstein metrics and smooth structures]{Einstein metrics 
and the number of smooth structures on a four--manifold}
\author{V.~Braungardt}
\address{Mathematisches Institut, Ludwig-Maximilians-Universit\"at M\"unchen,
Theresienstr.~39, 80333 M\"unchen, Germany}
\email{Volker.Braungardt@mathematik.uni-muenchen.de}
\author{D.~Kotschick}
\address{Mathematisches Institut, Ludwig-Maximilians-Universit\"at M\"unchen,
Theresienstr.~39, 80333 M\"unchen, Germany}
\email{dieter@member.ams.org}
\thanks{The second author is grateful to R.~Stern for pointing out 
the work of J.~Park. Support from the {\sl Deutsche Forschungsgemeinschaft}
is gratefully acknowledged. The authors are members of the 
{\sl European Differential Geometry Endeavour} (EDGE), Research 
Training Network HPRN-CT-2000-00101, supported by The European Human 
Potential Programme}
\date{February 2, 2003; MSC 2000 classification: primary 57R55; secondary 14J29, 
14J80, 53C25, 57R57}

\begin{document}

\begin{abstract}
We prove that for every natural number $k$ there are simply connected 
topological four--manifolds which have at least $k$ distinct smooth 
structures supporting Einstein metrics, and also have infinitely many 
distinct smooth structures not supporting Einstein metrics. Moreover, all 
these smooth structures become diffeomorphic to each other after connected 
sum with only one copy of the complex projective plane. We prove that 
manifolds with these properties cover a large geographical area.
\end{abstract}

\maketitle

\section{Introduction}

All the classical obstructions to the existence of Einstein metrics on 
four--manifolds are homotopy invariant. If a closed orientable 
$4$--manifold $M$ admits an Einstein metric, then its Euler characteristic 
has to be non--negative, and, furthermore, the Hitchin--Thorpe 
inequality
\begin{equation}\label{eq:HT}
e (M) \geq \frac{3}{2}\vert\sigma (M)\vert 
\end{equation}
must hold~\cite{HT}, where $e$ denotes the Euler characteristic and 
$\sigma$ the signature. This condition is clearly homotopy invariant, as 
are the restrictions coming from Gromov's notion of simplicial 
volume~\cite{gromov,K2}, and from the existence of maps of non--zero 
degree to hyperbolic manifolds~\cite{BCG}.

Using Seiberg--Witten invariants, LeBrun~\cite{lebrun} gave the first 
examples of simply connected smooth four-manifolds which satisfy the 
(strict) Hitchin--Thorpe inequality, but still do not admit Einstein 
metrics. As his examples were not known to be homeomorphic to manifolds 
admitting Einstein metrics, LeBrun's paper implicitly raised the question 
whether the new obstruction might in fact be homotopy invariant, or not. 
This issue was disposed of by the second author in~\cite{K}. Using 
LeBrun's~\cite{lebrun} work, Kotschick~\cite{K} proved the following result, 
showing for the first time that the smooth structures of $4$--manifolds form 
definite obstructions to the existence of an Einstein metric.
\begin{thm}\label{t:GT}
There are infinitely many pairs $(X_i,Z_i)$ of simply connected
closed oriented smooth $4$--manifolds such that:
\begin{itemize}
\item[{\bf 1)}] $X_i$ is homeomorphic to $Z_i$,
\item[{\bf 2)}] if $i\neq j$, then $X_i$ and $X_j$ are not homotopy
equivalent,
\item[{\bf 3)}] $Z_i$ admits an Einstein metric but $X_i$ does not,
\item[{\bf 4)}] $e (X_i) > \frac{3}{2}\vert\sigma (X_i)\vert$.
\end{itemize}
\end{thm}
\noindent
Note that 3)~implies in particular that $X_i$ and $Z_i$ are not diffeomorphic. 

After the proof of Theorem~\ref{t:GT}, Kotschick asked how many smooth 
structures with Einstein metrics and how many without such metrics exist 
on a given topological manifold, see~\cite{K} p.~6-7. He pointed out that, 
using for example the work of Fintushel--Stern, one can show that one has 
infinitely many choices for the smooth structures of the manifolds $X_i$ in 
Theorem~\ref{t:GT}. Kotschick~\cite{K} also remarked that by the work of 
Salvetti~\cite{Salve}, the number of distinct smooth structures among sets 
of homeomorphic minimal surfaces of general type can be arbitrarily large, 
and that all the examples in~\cite{Salve} have ample canonical bundle, and 
therefore have K\"ahler--Einstein metrics of negative scalar curvature. Thus, 
the number of smooth structures admitting Einstein metrics can be arbitrarily 
large.

The purpose of this paper is to show that these two phenomena, infinitely 
many smooth structures without Einstein metrics and an arbitrarily large 
number of smooth structures with Einstein metrics, can be realized on the 
same topological manifold. We shall prove the following:
\begin{thm}\label{t:main}
    For every natural number $k$ there is a simply connected topological 
    $4$-manifold $M_{k}$ which has at least $k$ distinct smooth structures 
    $Z_{k}^{i}$ supporting Einstein metrics, and also has infinitely many 
    distinct smooth structures $X_{k}^{j}$ not supporting Einstein metrics. 
    
    Moreover, all the $Z_{k}^{i}\#\C P^{2}$ and $X_{k}^{j}\#\C P^{2}$ for fixed 
    $k$ are diffeomorphic to each other.
    \end{thm}
We shall produce lots of such examples, with ratios $\vert\sigma\vert/e$ 
which are dense in the interval $[\tfrac{1}{3},\tfrac{1}{2}]$, compare 
Theorem~\ref{t:main2} in Section~\ref{s:proof}.

As the $Z_{k}^{i}$ and $X_{k}^{j}$ are all homeomorphic for fixed $k$, 
Wall's classical result~\cite{W} implies that they are stably diffeomorphic. 
That a single stabilization with $\C P^{2}$ suffices can be interpreted 
to mean that these differentiable structures are as close to each other as 
is possible while still being non-diffeomorphic. In fact, we shall exhibit 
$Z_{k}^{i}$ and $X_{k}^{j}$ as in Theorem~\ref{t:main} which are almost 
completely decomposable (ACD) in the sense of Mandelbaum~\cite{Msurv}, so 
that their connected sums with $\C P^{2}$ are diffeomorphic to 
$p\C P^{2}\# q\overline{\C P^{2}}$ for some $p$ and $q$. Whether such a 
decomposable manifold can admit an Einstein metric is only known in very 
few cases with $p=1$.

To put Theorem~\ref{t:main} into perspective, we continue with the 
chronology of earlier work in this direction. After~\cite{K} appeared, 
LeBrun~\cite{lebrun2,lebrun3} refined his arguments from~\cite{lebrun}, 
and produced more examples of precisely the type exhibited in 
Theorem~\ref{t:GT}, where one has pairs of homeomorphic manifolds such that 
one is Einstein and the other is not. However, he did not discuss the 
number of smooth structures. This was taken up recently by Ishida and LeBrun 
in~\cite{IL}. They give examples of simply connected topological four-manifolds
with an infinite number of smooth structures which do not admit Einstein 
metrics. Like LeBrun in his earlier papers~\cite{lebrun,lebrun2,lebrun3}, 
Ishida--LeBrun~\cite{IL} do not exhibit multiple smooth structures with 
Einstein metrics on the same manifold where one has infinitely many smooth 
structures without such metrics. In fact, for the most interesting ones of 
their examples, no smooth structure with an Einstein metric is known.

One of the difficulties in proving results like Theorems~\ref{t:GT} 
and~\ref{t:main} above is that there are almost no existence results 
for Einstein metrics on simply connected $4$-manifolds. Therefore, one 
is always forced to arrange a situation where one can appeal to the only 
existence result covering lots of homeomorphism types, which is the 
resolution of the Calabi conjecture for negative scalar curvature due to 
Aubin~\cite{A} and Yau~\cite{Y}. This then leads to questions about the 
geography of complex surfaces of general type, and of some related classes 
of four-manifolds. Thus, in the present paper we make substantial progress 
on two geographical questions, which are of interest independently of the 
applications to Einstein metrics. One is the geography of algebraic
surfaces which are iterated branched covers of the plane, the other is 
about symplectic four-manifolds which are almost completely decomposable. 

Salvetti~\cite{Salve} considered iterated cyclic branched covers of the 
projective plane, and used these to prove that for any $k$, there exists 
a pair of invariants $e$ and $\sigma$ such that for this pair one has at 
least $k$ homeomorphic surfaces with different divisibilities for their 
canonical classes. In his examples, the ratios $\sigma/e$ are so close 
to zero that one cannot use them to prove Theorem~\ref{t:main} with 
the arguments of~\cite{lebrun,K}. In fact, even the improved estimates 
of~\cite{lebrun2,lebrun3} do not apply. Therefore, in Section~\ref{s:salve} 
below we provide a generalization of Salvetti's arguments which 
shows that, by choosing the parameters judiciously, iterated cyclic 
branched covers of the projective plane can be used to cover other 
parts of the geography of surfaces. In particular, we can arrange 
$k$--tuples of homeomorphic surfaces with different divisibilities for 
the canonical class with characteristic numbers which are such that 
homeomorphic manifolds without Einstein metrics can be found using the 
improved estimate from~\cite{lebrun3}.

In Section~\ref{s:Park}, which is inspired in part by the work of 
J.~Park~\cite{P2}, we discuss the geography of minimal symplectic 
four-manifolds which are almost completely decomposable. Blowups of 
these will be used for the $X_{k}^{j}$ in Theorem~\ref{t:main}. Even 
without the ACD requirement, our geography results are qualitatively 
stronger than what was known before, compare for example~\cite{Go,P2,P3}.

In Section~\ref{s:proof} we combine the different ingredients to prove a 
more precise version of Theorem~\ref{t:main}. We shall also exhibit 
infinitely many smooth structures without Einstein metrics on many other 
manifolds which are not homotopy equivalent to complex surfaces, for which 
the existence of smooth structures with Einstein metrics is an open question. 
See Theorem~\ref{t:pq}.

In Section~\ref{s:ex} we give explicit examples of manifolds with very 
small homology which have a smooth structure supporting an Einstein 
metric and have infinitely many smooth structures which do not support such a 
metric. We also give simple explicit examples with multiple smooth structures 
admitting Einstein metrics. 

\section{The geography of iterated branched covers of $\C P^{2}$}\label{s:salve}

In this section we study the spread of Chern numbers among algebraic 
surfaces which are iterated branched covers of the projective plane.
We build on the work of Salvetti~\cite{Salve} to show that for any 
integer $k$ there are $k$-tuples of homeomorphic surfaces with ample 
canonical classes of different divisibilities, whose ratio $c_1^2/\chi$ 
can be specified arbitrarily within a certain range. Note that 
$k$-tuples of homeomorphic surfaces with canonical classes of different 
divisibilities were first exhibited by Catanese~\cite{Ca2} using bidouble 
covers of $\C P^{1}\times\C P^{1}$. In his examples the divisibilities are 
even, but the method probably extends to odd divisibility. Nevertheless, 
we found iterated covers of the plane to be more convenient to use.

Given positive integers $r$, $d_1,\ldots,d_r$ and $m_1,\ldots,m_r$,
one can construct a simply-connected complex algebraic surface $S$ 
by starting from the projective plane and repeatedly passing to coverings 
of degrees $d_j$ branched along the preimages of smooth curves of degree 
$n_j=d_jm_j$ in the plane. The canonical class of $S$ is 
$\sum_{j=1}^r(d_j-1)m_j-3$ times the pullback of the class of a line in 
the projective plane. Except for some small values of the parameters, 
the surface $S$ so obtained is minimal of general type and has ample 
canonical bundle. The Chern numbers of $S$ are
\begin{eqnarray*}
c_1^2(S)&=&d_1\cdot\ldots\cdot d_r(\sum_{j=1}^r(d_j-1)m_j-3)^2 \ 
,\label{c12}\\
c_2(S)&=&\frac12d_1\cdot\ldots\cdot d_r\big((\sum_{j=1}^r(d_j-1)m_j-3)^2+
 (\sum_{j=1}^r(d_j^2-1)m_j^2-3)\big) \ ,
\end{eqnarray*}
and they determine the holomorphic Euler characteristic $\chi(S)$ and 
the signature
\begin{equation*}
\sigma(S)=-\frac13d_1\cdot\ldots\cdot d_r(\sum_{j=1}^r(d_j^2-1)m_j^2-3) \ .
\end{equation*}

We consider the inverse problem, starting with a fixed pair of 
invariants, say $c_1^2(S)$ and $\sigma(S)$, and try to find $k$ solutions 
of the above equations for $d_1,\ldots,d_r$ and $m_1,\ldots,m_r$. 
Salvetti~\cite{Salve} considered the special case when the covering degrees 
$d_j$ are all equal and the $m_j$ are not too far from being equal; more 
precisely he assumed $\sum_{j=1}^rm_j^2\le\tfrac1{r-1}(\sum_{j=1}^rm_j)^2$. 
This leads to ratios for $c_1^2(S)/\chi(S)$ close to $8$, which is not 
suitable for our purposes. In order to use these surfaces as the $Z_{k}^{i}$ 
in Theorem~\ref{t:main}, we need $c_1^2(S)/\chi(S)$ to be somewhere below $6$, 
and the smaller we get this ratio, the easier the proof will be.
In order to minimize $c_1^2(S)/\chi(S)$ we have to maximize the quotient of
$\sum(d_j^2-1)m_j^2$ by $(\sum(d_j-1)m_j)^2$, i.\,e.\ we should have a few
of the $m_j$ much bigger than the others and the corresponding covering 
degrees $d_j$ small.

We can do better and adjust $c_1^2(S)/\chi(S)$ to approximate any value between
$4$ and $8$. 
\begin{sa}\label{ratio-4.8}
Let $k$ be a positive integer. There are values for $c_1^2$ and $\sigma$ which 
are realized by at least $k$ iterated branched covers with different 
divisibilities of the canonical class. The divisibilities can be arranged to 
be all even, or all odd. The corresponding ratios $c_1^2/\chi$ are dense 
in the interval $[4,8]$.
\end{sa}
\begin{proof}
To spread out the Chern numbers, fix rational numbers $\mu_1,\ldots,\mu_s$ 
normalized by $\mu_1+\ldots+\mu_s=1$. Put $\mu^2=\sum_{j=1}^s\mu_j^2$ and
note that by a suitable choice of $s$ and $\mu_j$ we can place $\mu^2$ anywhere
between $0$ and $1$ because the set of numbers of the form
$\sum_{j=1}^s\mu_j^2$ with arbitrary $s$ and rational
$\mu_1,\ldots,\mu_s$ summing up to 1 is dense in the unit interval.

Consider a tower of coverings of the projective plane by $s$ iterated double 
covers branched over curves of degree $2\mu_jm_0$, where the integer $m_0$ will 
be fixed later. Of course $m_{0}$ has to be a multiple of the denominators of
$\mu_1,\ldots,\mu_s$. In the end we will let $m_0$ grow to infinity.

Now on top of this tower we consider $16$ further cyclic covers of very
high degree $d$ branched over the preimages of curves of degrees 
$dm_1,\ldots,dm_{16}$. If $m_0$ is sufficiently large, the elementary number 
theory worked out by Salvetti~\cite{Salve} will provide us with
several solutions for  $(d,m_1,\ldots,m_{16})$ giving rise to the same 
invariants $c_1^2$ and $\sigma$ for the total covering surfaces.

For our $s+16$ stage tower, the formulae for the invariants specialize to
\begin{eqnarray}
c_1^2(S)&=&2^{s}d^{16}(m_0+(d-1)\sum_{j=1}^{16}m_j-3)^2\label{s}\\
\sigma(S)&=&-\frac{1}{3}2^{s}d^{16}(3\mu^2m_0^2+
(d^2-1)\sum_{j=1}^{16}m_j^2-3) \ .\label{ss}
\end{eqnarray}

Fix $\delta>0$. Since the sum $\sum_pp^{-1}$ over prime numbers $p$ diverges, 
for any real number $\alpha>1$ the number of primes between $\alpha^n$ and 
$\alpha^{n+1}$ will be unbounded for $n\to\infty$. Hence we can find $k+1$ odd 
primes such that the largest one is at most $\alpha$ times the smallest one. 
Forgetting the smallest one, we obtain a set $D$ of $k$ odd primes such that 
$(d+1)/(d'-1)<\alpha$ for any $d,d'\in D$. Note that $D$ depends on $\alpha$, 
though this is not explicit in our notation. Put $d^\ast=\min(D)$ and 
$\eps=\delta/(d^\ast-1)$. It is clear that $\eps$ converges to $0$ for 
$\alpha\to1$ and $\delta$ fixed.

These $d\in D$ will be used as degrees in our tower of coverings. Put 
$P=\prod_{d\in D}d$. By~\eqref{s} and since the $d\in D$ are primes, we can 
write $c_1^2=2^{s}P^{16}Q$ for some integer $Q$. Furthermore, since 
$c_1^2/2^{s}d^{16}$ has to be a square, we can write $Q=C^2$ for an integer
$C$.  Similarly, we can write $\sigma=-\tfrac132^{s}P^{16}C'$ for some integer
$C'$.  The equations~\eqref{s} and~\eqref{ss} now read 
\begin{eqnarray}
(d-1)\sum_{j=1}^{16}m_j&=&(P/d)^8C-(m_0-3) \ ,\label{sss}\\
(d^2-1)\sum_{j=1}^{16}m_j^2&=&(P/d)^{16}C'-3(\mu^2m_0^2-1) \ .\label{ssss}
\end{eqnarray}
We are left with the task of solving the pair of equations
\begin{align}\label{to-be-solved}
\sum_{j=1}^{16}m_j&=A_d,&\sum_{j=1}^{16}m_j^2&=B_d
\end{align}
for each $d\in D$ (separately), where we have put
\begin{align*}
A_d&=((P/d)^8C-(m_0-3))/(d-1) \ ,\\
B_d&=((P/d)^{16}C'-3(\mu^2m_0^2-1))/(d^2-1) \ .
\end{align*}
We can achieve that both $A_d$ and $B_d$ are integers by the following
\begin{hs}\label{C-and-C'-exist}
There are integers $C,C'$ such that for every $d\in D$ we have
$(P/d)^8C\equiv m_0-3$ {\rm mod} $d-1$
and $(P/d)^{16}C'\equiv3(\mu^2m_0^2-1)$ {\rm mod} $d^2-1$.
\end{hs}
\begin{proof}
Since $D$ consists of nearby primes we infer that $d$ and $d'\pm1$ are 
coprime for any $d,d'\in D$.  Hence $P$ is invertible $\mod\prod_d(d-1)$ 
and we can find $C$ with $P^8C\equiv m_0-3$ mod $\prod_d(d-1)$.
Since $(P/d)^8\equiv P^8$ mod $d-1$ for every $d$ the result follows.

Similarly $P$ is invertible mod $\prod_d(d^2-1)$ so we find $C'$
with $P^{16}C'\equiv3(\mu^2m_0^2-1)$ mod $\prod_d(d^2-1)$.
Since $d^{16}\equiv 1$ mod $d^2-1$ the result follows.
\end{proof}
We are free to modify $C$ by a multiple of $\prod_d(d-1)$ and $C'$
by a multiple of $\prod_d(d^2-1)$; hence we can arrange that $A_{d^\ast}$
and $B_{d^\ast}$ have distance less than $P^9$ respectively $P^{18}$
from any given values $A$ and $B$.  Notice that with $C$ and $C'$
satisfying Lemma~\ref{C-and-C'-exist} both $A_d$ and $B_d$ will be even
for every $d$ (as soon as $k\ge2$).

Our goal is to choose $C,C'$ in such a way that we can solve
the pair of equations \eqref{to-be-solved} relying on the following
\begin{hs}[Salvetti~\cite{Salve}]\label{salvetti-lemma}
Let $A,r$ be integers with $A>0$ and $r\ge16$.
The integral quadratic form $x_1^2+\ldots+x_r^2$, under the restriction
$x_1+\ldots+x_r=A$ and $x_j>0$ for $j=1,\ldots,r$,
represents all integers $B\equiv A$ $(\text{\rm mod}\;2)$ such that
$$A^2/r+\alpha_r\le B\le A^2/(r-1),$$
where $\alpha_r$ is a constant depending only on $r$.
\end{hs}
\begin{proof}
This is the Lemma on page 166 of~\cite{Salve}. The proof is elementary, 
using that each non-negative integer is a sum of four squares.
\end{proof}
Recall that we denote the smallest element of $D$ by $d^\ast$.
For any positive integer $m_0$ choose $C(m_0)$ such that 
$$
A_{d^\ast}=((P/d^\ast)^8C(m_0)-(m_0-3))/(d^\ast-1)
$$
satisfies $|A_{d^\ast}-\eps m_0|<P^9$.  Since $P$ does not depend on $m_0$, 
we will eventually have
$$
\tfrac12\eps^2m_0^2+240\alpha_{16}<A_{d^\ast}^2<(2\eps m_0)^2
$$
for all sufficiently large $m_0$, where $\alpha_{16}$ is the constant 
from Lemma~\ref{salvetti-lemma} with $r=16$. Then the quantity
$$
\Delta(m_0)=\frac{A_{d^\ast}^2}{15}-(\frac{A_{d^\ast}^2}{16}+\alpha_{16})
  =\frac{A_{d^\ast}^2}{240}-\alpha_{16}
$$
will be bounded below by
\begin{equation}\label{Delta-lower-bound}
\Delta(m_0)>\frac{\eps^2}{480}m_0^2
\end{equation}
for large $m_0$.
Now choose $C'(m_0)$ such that $B_{d^\ast}$ differs from
$A_{d^\ast}^2/16+\Delta(m_0)/2$ by no more than $P^{18}$.
Then for large $m_0$ we will have
\begin{equation}\label{salvetti-criterion}
\frac{A_{d^\ast}^2}{16}+\alpha_{16}+\frac13\Delta(m_0)<
B_{d^\ast}<\frac{A_{d^\ast}^2}{15}-\frac13\Delta(m_0) \ .
\end{equation}

\begin{hs}
If $\alpha$ is sufficiently close to 1 then for all $m_0\gg0$ we have
$$\frac{A_d^2}{16}+\alpha_{16}<B_d<\frac{A_d^2}{15}$$
for every $d\in D$.
\end{hs}
\begin{proof}
We first show that $A_d^2$ differs from $A_{d^\ast}^2$ by no more
than $\Delta(m_0)/6$.  To see this, observe that
\begin{alignat*}{1}
|(d-1)A_d-(d^\ast-1)A_{d^\ast}|
  &=(P/d^\ast)^8C(m_0)|(d^\ast/d)^8-1|\\
  &\le(\alpha^8-1)((d^\ast-1)A_{d^\ast}+m_0-3) \ .
  \end{alignat*}
We divide by $d-1$ and use $|(d^\ast-1)/(d-1)-1|<\alpha-1$ to obtain
$$
|A_d-A_{d^\ast}|\le(\alpha^8-1)(A_{d^\ast}+
\frac{m_0-3}{d-1})+(\alpha-1)A_{d^\ast}$$
or, using $A_{d^\ast}<2\eps m_0$ and $1/(d-1)\le\eps/\delta$,
$$|A_d-A_{d^\ast}|\le h(\alpha)\eps m_0$$
with $h(\alpha)=(\alpha^8-1)(2+1/\delta)+2(\alpha-1)$, which goes
to zero when $\alpha\to1$.  Then
$$|A_d^2-A_{d^\ast}^2|\le2A_{d^\ast}|A_d-A_{d^\ast}|+|A_d-A_{d^\ast}|^2\le 
h'(\alpha)\eps^2m_0^2$$
where $h'(\alpha)=4h(\alpha)+h(\alpha)^2$ also goes to zero when $\alpha\to1$.
This shows that as soon as $\alpha$ is closer to 1 than some constant depending
only on $\delta$, the difference $|A_d^2-A_{d^\ast}^2|$ is less than
1/6 times the term on the right-hand side of~\eqref{Delta-lower-bound}.

A similar computation shows that we can assume the same bound for
$|B_d-B_{d^\ast}|$.  In detail,
\begin{alignat*}{1}
|(d^2-1)B_d-({d^\ast}^2-1)B_{d^\ast}|
  &\le(P/d^\ast)^{16}C'(m_0)|(d^\ast/d)^{16}-1|\\
  &\le(\alpha^{16}-1)(({d^\ast}^2-1)B_{d^\ast}+3\mu^2m_0^2-3) \ .
  \end{alignat*}
Dividing by $d^2-1$ we obtain
$$|B_d-B_{d^\ast}|
  \le(\alpha^{16}-1)(B_{d^\ast}+\frac{3\mu^2m_0^2-3}{d^2-1})+(\alpha^2-1)B_{d^\ast} 
  \ .
  $$
Using $B_{d^\ast}<\frac1{15}A_{d^\ast}^2<\frac4{15}\eps^2m_0^2$ and
$1/(d^2-1)\le(\eps/\delta)^2$ this gives
$|B_d-B_{d^\ast}|<h''(\alpha)\eps^2m_0^2$ where
$h''(\alpha)=(\alpha^8-1)(\frac4{15}+3\mu^2/\delta^2)+\frac4{15}(\alpha^2-1)$ 
becomes arbitrarily small as $\alpha\to1$.

Now the claim follows from the estimates~\eqref{salvetti-criterion} because
replacing $d^\ast$ by $d$ does not change any of the terms by more than
$\Delta(m_0)/6$.
\end{proof}

According to Lemma~\ref{salvetti-lemma} this shows that for each $d\in D$
the pair of equations~\eqref{to-be-solved} is solvable by positive integers. 

Next we calculate the ratios $c_1^2/\chi$ for these surfaces. As they 
all have the same invariants, we can look at the one with the reference parameter 
$d=d^\ast$. For the degrees of the branch divisors we have the estimate 
$\sum_{j=1}^rm_j=A_d\le2\eps m_0$.
This gives
\begin{alignat*}{1}
c_1^2 &=2^{s}d^{16}(m_0+(d-1)\sum_{j=1}^{16}m_j-3)^2 \\
&=2^{s}d^{16}m_0^2(1+(d-1)O(\eps))^2
\;=\;2^{s}d^{16}m_0^2(1+O(\delta))^2 \ ,
\end{alignat*}
and, using $\sum_jm_j^2\le(\sum_jm_j)^2$,
\begin{eqnarray*}
-\sigma &=&\tfrac132^{s}d^{16}(3\mu^2m_0^2+(d^2-1)\sum_{j=1}^{16}m_j^2-3)\\
&=&2^{s}d^{16}\mu^2m_0^2(1+(d^2-1)O(\eps^2))
\;=\;2^{s}d^{16}\mu^2m_0^2(1+O(\delta^2)) \ .
\end{eqnarray*}
Since we can choose $\delta$ arbitrarily small, then $\alpha$ arbitrarily
close to $1$ (increasing $d$) and, finally, $m_0$ arbitrarily large, 
these estimates for $c_1^2$ and $\sigma$ show that we can arrange 
$-\sigma/c_1^2$ arbitrarily close to $\mu^2$.

Recall that the possible values of $\mu^2$ are dense in the unit interval. 
Thus, $\frac{c_{1}^{2}}{\chi}=\frac{8}{1-\sigma/c_{1}^{2}}$ ranges over 
values dense in $[4,8]$. 

It remains to check the divisibilities of the canonical classes of the 
surfaces with fixed invariants that we have constructed. The divisibility 
is $m_0+(d-1)\sum_jm_j-3$, because the pullback of the hyperplane 
class to the iterated covering is primitive, compare~\cite{Mjdg}
Proposition 10 and corollary. Thus 
the divisibility is odd if and only if $m_{0}$ is 
even. This imposes no restriction on our construction. We obtain even 
divisibilities if and only if $m_0$ is odd. Since $m_0$ is restricted to 
multiples of the denominators of the $\mu_j$, this requires that these 
denominators be odd. This restriction still leaves us with a set of attainable 
$\mu^2$ which is dense in the unit interval.

Now for each $d\in D$ we obtain a unique total covering degree $2^{s}d^{16}$.
This implies that the surfaces obtained for a fixed $c_{1}^{2}$ and different
$d$ have different divisibilities $d^8\sqrt{c_1^2/2^{s}}$ for their canonical
classes.
\end{proof}

\begin{cor}\label{c:geography}
    For every $k>0$ there are $k$-tuples of simply connected spin and 
    non-spin complex algebraic surfaces with ample canonical bundles which are 
    homeomorphic, but are pairwise non-diffeomorphic. For every $k$, 
    the ratios $c_1^2/\chi$ of such $k$-tuples are dense in the interval $[4,8]$.
\end{cor}
\begin{proof}
    The surfaces in Theorem~\ref{ratio-4.8} are all simply connected and have 
    ample canonical bundles. They are spin or non-spin according to whether 
    the divisibility of the canonical class is even or odd. Once the parity of 
    the divisibility and the values of the Chern numbers are fixed, all these 
    surfaces are homeomorphic by Freedman's result~\cite{freed}. However, 
    surfaces with different divisibilities cannot be diffeomorphic, because 
    Seiberg--Witten theory shows that the canonical class of a minimal surface 
    of general type is diffeomorphism-invariant (up to sign), 
compare~\cite{Wi}, page 789.
\end{proof}

\section{Geography of ACD symplectic manifolds}\label{s:Park}

In this section we study the geography of simply connected minimal 
symplectic four--manifolds which are almost completely decomposable. 
Recall that a four--manifold $X$ is called almost completely decomposable 
or ACD if $X\#\C P^{2}$ is diffeomorphic to $p\C P^{2}\# q\overline{\C P^{2}}$ 
for some $p$ and $q$. Mandelbaum~\cite{Msurv} conjectured that every 
simply connected complex algebraic surface is ACD, and it is very 
natural to extend this to simply connected symplectic four--manifolds. 
Mandelbaum and Moishezon proved that certain algebraic surfaces, 
including simply connected elliptic surfaces, complete intersections 
and double planes are indeed ACD, compare~\cite{Ma,Msurv,MM} and the 
references cited there.

The geography of minimal symplectic four--manifolds has been investigated 
by many authors in recent years, for example by Fintushel--Stern~\cite{FScusp}, 
Gompf~\cite{Go1}, Park~\cite{P2,P3} and Stipsicz~\cite{BMY}. Here we reprove and 
improve their results, although we only use manifolds with the ACD property. 
That the ACD condition does not constrain the geography can be taken as evidence 
that all simply connected symplectic manifolds may be ACD. The restriction to 
minimal symplectic manifolds is natural, and implies that all our manifolds 
will be irreducible, see~\cite{Bourbaki}, which gives the ACD 
property added interest.

We use the coordinates $c_{1}^{2}=2e+3\sigma$ and $\chi=\frac{1}{4}(e+\sigma )$ 
to state our geography results. By the work of Taubes~\cite{Taubes,Bourbaki}, 
minimal symplectic four--manifolds satisfy $c_{1}^{2}\geq 0$. It is clear that 
in the simply connected case one must have $\chi > 0$. Thus, we try to cover 
lattice points in the first quadrant of the $(\chi , c_{1}^{2})$--plane with
simply connected minimal symplectic manifolds which dissolve after connected
sum with only one copy of $\C P^{2}$. All our manifolds are, or can be made,
non-spin, and we will not repeat this. We shall use symplectic summation 
along submanifolds, as pioneered by Gompf~\cite{Go1}, with only a handful of 
building blocks. Most of these summations will be along tori of zero 
selfintersection, in which case the Chern invariants are additive. For the 
resulting manifolds we have to check minimality and the ACD property. 
Minimality is always true if we sum minimal symplectic manifolds with 
$\chi >1$. Proofs of this have been given by  Li--Stipsicz~\cite{LS} and by 
Park~\cite{P2}, with the latter attributing the result to Lorek. To prove 
almost complete decomposability we shall use the following ``irrational 
connected sum lemma'' of Mandelbaum, see~\cite{Ma,Msurv,Go}.
\begin{prop}\label{p:diss}
    Let $M$ and $N$ be simply connected oriented $4$-manifolds
    containing the same embedded surface $F$ of genus $g\geq 1$ with zero
    selfintersection. Assume that $F$ has simply connected complement
    in $M$. Denote by $P$ the sum of $M$ and $N$ along $F$, and assume 
    that $P$ is not spin. 
    \begin{enumerate}
    \item Then $P\# \C P^{2}\# \overline{\C P^{2}}$ is diffeomorphic
    to $M\# N\# 2g(\C P^{2}\# \overline{\C P^{2}})$.
    \item If $(N,F)$ is obtained from a pair $(N',F')$ by blowing up a 
    point on $F'$, then $P\# \C P^{2}$ is diffeomorphic
    to $M\# N'\# 2g(\C P^{2}\# \overline{\C P^{2}})$.
    \end{enumerate}
    \end{prop}

We now list our building blocks.
\begin{ex}[Elliptic building blocks]
We shall denote by $E(n)$ the relatively minimal elliptic surface with 
$\chi (E(n))=n$ over $S^{2}$ without multiple fibers. This is the fiber 
sum of $n$ copies of $E(1)=\C P^{2}\# 9\overline{\C P^{2}}$, and is 
therefore ACD by the second part of Proposition~\ref{p:diss}. It is minimal 
as soon as $n>1$.

We can use logarithmic transformations on $E(n)$ to produce infinitely 
many distinct smooth structures, all of which support symplectic structures 
such that the fibers are symplectic submanifolds. Using such logarithmic 
transformations we can also change the homeomorphism type of $E(2n)$, which 
is spin, to a non-spin elliptic surface. All these elliptic surfaces with 
multiple fibers are ACD, see~\cite{Ma,MM}.

Gompf~\cite{Go1} has shown that the $K3$ surface $E(2)$ contains two 
disjoint nuclei corresponding to different elliptic fibrations in such 
a way that there is a symplectic form for which the tori in the two 
nuclei are simultaneously symplectic. This is useful because each of 
the nuclei contains a $2$-sphere intersecting the torus fiber once, so 
that one can perform logarithmic transformations or symplectic 
summations independently inside the two nuclei, without introducing 
a nontrivial fundamental group.
    \end{ex}
\begin{ex}[A small building block]\label{bb:s}
We shall denote by $S$ the simply connected symplectic manifold $S_{1,1}$ 
constructed by Gompf in~\cite{Go1}, Example 5.4. It has $c_{1}^{2}(S)=1$, and 
$\chi (S)=2$, and contains a symplectically embedded torus $T$ of zero 
selfintersection and a symplectically embedded genus $2$ surface $F$
disjoint from $T$, such that $S\setminus (T\cup F)$ is simply connected. 
That $S$ is irreducible was proved by Stipsicz~\cite{SS}. {\it A fortiori} 
it is minimal.
    \end{ex}
\begin{ex}[Building blocks with positive signature]\label{bb:pos}
We use the following construction of Li--Stipsicz~\cite{LS}, compare 
also~\cite{BMY}. 

For every positive integer $n$ there is a symplectic manifold $X_{n}$ 
which is a Lefschetz fibration over the surface $\Sigma_{n+2}$ of 
genus $n+2$ which admits a section of selfintersection $-n-1$. 
It has Chern invariants
$$
\chi (X_{n})=25n^{2}+30n+1  \ ,
$$
$$
c_{1}^{2}(X_{n})=225n^{2}+180n \ .
$$
Furthermore, the fibration and the section induce inverses of each other 
on the fundamental groups. Thus one can kill the fundamental group 
of $X_{n}$ by symplectic summation along the section. We will use a 
blowup of the $K3$ surface as follows. First we construct a smooth 
symplectic submanifold inside a nucleus by smoothing the union of $n+2$ 
copies of a regular fiber and one copy of a section. This gives a 
surface of genus $n+2$ and selfintersection $2n+2$. Blowing up $n+1$ 
points on this surface, its selfintersection drops to $n+1$, so that 
it can be symplectically summed to the section of $X_{n}$. Note that 
the surface has simply connected complement inside the blown-up 
nucleus of $K3$. In this way we obtain a simply connected minimal 
symplectic $4$--manifold $Y_{n}$. Using the above formulae for the Chern 
invariants of $X_{n}$ we obtain:
$$
\chi(Y_{n})=25n^{2}+31n+4 \ ,
$$
$$
c_{1}^{2}(Y_{n})=225n^{2}+187n+7 \ .
$$
These manifolds are not spin and still contain a nucleus of a $K3$ surface with 
simply connected complement. Note that for every $\epsilon>0$ there 
is an $n$ such that $c_{1}^{2}(Y_{n})/\chi(Y_{n})>9-\epsilon$.
\end{ex}

As a warmup for our main geography result we first show how to fill up a 
certain region, which includes that below the Noether line. The importance 
of this is that the width of this region in the $y$ direction goes to infinity 
with $x$.

For a constant $c$ let $R_{c}$ denote the set of lattice points 
$(x,y)$ in the plane satisfying $x>0$, $y\geq 0$, and 
\begin{alignat}{1}
    y &\leq 3x-51 \ , \\
    y &\leq 6x-c \ . \label{c}
    \end{alignat}

\begin{prop}\label{p:wedge}
There exists a constant $c$ such that all lattice points in $R_{c}$ 
are realized as the Chern invariants $(\chi,c_{1}^{2})$ of infinitely 
many homeomorphic pairwise nondiffeomorphic simply connected minimal 
symplectic manifolds, all of which are almost completely decomposable.	
\end{prop}
\begin{proof}
    All our examples will be non-spin and will have the same Chern 
    invariants. Thus they are homeomorphic by Freedman's 
    classification~\cite{freed}.
    
    There is nothing to prove for $y=0$, as minimal elliptic surfaces and 
    their logarithmic transformations give the required examples. Thus we 
    may assume $y>0$. Given a positive integer $y$, we can write it uniquely 
    as $y=9k+r-8$, with $0\leq r\leq 8$ and $k>0$. Then consider the manifold 
    $X(k,r,n)$ obtained as the symplectic sum of $k$ copies of building block 
    $S$ summed along the genus $2$ surface $F$, and of $r$ further copies of 
    $S$ and one copy of $E(n)$ summed to the result along the torus $T$. This 
    is again simply connected. The Chern invariants are 
    $c_{1}^{2}(X(k,r,n))=9k+r-8$ and $\chi(X(k,r,n))=3k+2r+n-1$. If we take 
    $n\geq 2$, then the building blocks are minimal, and so are the $X(k,r,n)$. 
    Moreover, the $X(k,r,n)$ fill out the claimed region (for any $c$).
    
    Consider now the connected sum $X(k,r,n)\#\C P^{2}$. Applying the 
    second part of Proposition~\ref{p:diss} to the seam inside the 
    elliptic piece $E(n)=E(1)\cup_{T^{2}} E(n-1)$ we can split off a copy of 
    $\C P^{2}\# \overline{\C P^{2}}$. Then using this to apply the 
    first part of Proposition~\ref{p:diss} to the remaining seams, 
    and breaking up the elliptic pieces, we see that 
    \begin{alignat*}{1}
    X(k,r,n) &\#\C P^{2} \\
    &\cong (k+r)S\# (3k+r+2n-2)\C P^{2}\# (3k+r+10n-3)
    \overline{\C P^{2}} \\
    &\cong (k+r)S\# (3k+r+2n-2)(S^{2}\times S^{2})\# (8n-1)\overline{\C P^{2}}  \ .
    \end{alignat*}
     By the result of Wall~\cite{W} there is a $k_{0}$ such that 
    $S\# k_{0}(S^{2}\times S^{2})$ is completely decomposable. Therefore, 
    $X(k,r,n)\# \C P^{2}$ dissolves as soon as $3k+r+2n-2\geq 
    k_{0}$, which follows from~\eqref{c} with $c=3k_{0}+72$.
  
    It remains to show that there are infinitely many symplectic
    manifolds homeomorphic but non-diffeomorphic to $X(k,r,n)$, all of
    which are ACD. For this we replace the elliptic surface $E(n)$ without 
    multiple fibers by one with multiple fibers obtained by logarithmic
    transformation. In this case the general fiber becomes 
    divisible in homology, in particular its complement is no longer simply
    connected. Here this is irrelevant because the torus in $S$ has 
    simply connected complement, so that the symplectic sum does give 
    a simply connected manifold and Proposition~\ref{p:diss} can be 
    applied.
    
    The logarithmic transformations on $E(n)$ produce infinitely many
    distinct smooth structures on the topological manifold underlying
    $E(n)$, which are detected by Seiberg--Witten invariants,
    cf.~\cite{Wi}. This difference in the Seiberg--Witten invariants
    survives the symplectic sum operation along a fiber, because of the
    gluing formulas due to Morgan--Mr\'owka--Szab\'o~\cite{MMS} and
    Morgan--Szab\'o--Taubes~\cite{MST}. Thus, we can produce infinitely
    many minimal symplectic manifolds homeomorphic but non-diffeomorphic
    to $X(k,r,n)$. All these are ACD by the same argument as for
    $X(k,r,n)$ (and the fact that the elliptic building blocks are ACD 
    even when they contain multiple fibers).  
\end{proof}

\begin{thm}\label{t:geo}
    For every $\epsilon >0$, there is a constant $c(\epsilon)>0$ such 
    that every lattice point $(x,y)$ in the first quadrant satisfying
    \begin{equation}\label{cutoff}
    y\leq (9-\epsilon)x-c(\epsilon)
    \end{equation}
    is realized by the Chern invariants $(\chi,c_{1}^{2})$ 
    of infinitely many pairwise non-diffeomorphic simply connected minimal 
    symplectic manifolds, all of which are almost completely decomposable.
    \end{thm}
    \begin{proof}
	Given $\epsilon>0$, we first choose an $i$ such that 
	$c_{1}^{2}(Y_{i})/\chi(Y_{i})>9-\epsilon$ for the building block 
	$Y_{i}$ in Example~\ref{bb:pos}. Denote this fixed $Y_{i}$ by $Y$.
	
        Let $Y(l,k,r,n)$ be the symplectic manifold obtained by symplectically 
	summing $l$ copies of $Y$ along the torus $T$ in the $K3$ nucleus 
	inside $Y$, and then summing the result to the manifold $X(k,r,n)$ from 
	the proof of Proposition~\ref{p:wedge} using the same torus in $Y$ and 
	the torus in $X(k,r,n)$ coming from the elliptic piece. If we choose 
	$c(\epsilon)$ large enough, then all the lattice points 
	satisfying~\eqref{cutoff} are covered by the translates of $R_{c}$ which 
	we obtain in this way. In all these summations the complement of the 
	surface along which the summation is performed is simply connected in at 
	least one of the summands, so that the resulting manifolds are simply 
	connected. They are all minimal, as we may assume $n>1$. 
        
        Consider now the connected sum $Y(l,k,r,n)\#\C P^{2}$. Applying the
	second part of Proposition~\ref{p:diss} to the seam inside the elliptic 
	piece $E(n)=E(1)\cup_{T^{2}} E(n-1)$ we can split off a copy of 
	$\C P^{2}\# \overline{\C P^{2}}$. Then using this to apply the first 
	part of Proposition~\ref{p:diss} to the remaining seams, and breaking 
	up the elliptic pieces, we see that 
\begin{alignat*}{1}
&Y(l,k,r,n) \#\C P^{2} \\
&\cong lY\# (k+r)S\# (l+3k+r+2n-2)\C P^{2}\# (l+3k+r+10n-3) \overline{\C P^{2}} \\
&\cong lY\# (k+r)S\# (l+3k+r+2n-2)(S^{2}\times S^{2})\# (8n-1)\overline{\C P^{2}}  \ .
\end{alignat*}
        By choosing $c(\epsilon)$ large enough, we can ensure that 
        $l+3k+r+2n-2$ is always larger than the ``resolving number'' 
        of $Y$ and of $S$, cf.~\cite{Msurv}. This means that the result of 
        Wall~\cite{W} can be applied to show that the above connected 
        sum is completely decomposable. Thus, if the Chern invariants 
	$(x,y)$ of $Y(l,k,r,n)$ satisfy $y\leq (9-\epsilon)x-c(\epsilon)$ 
	with large enough $c(\epsilon)$, we conclude that $Y(l,k,r,n)$ is ACD.
	
	It remains to show that there are infinitely many symplectic manifolds 
	homeomorphic but non-diffeomorphic to $Y(l,k,r,n)$, all of which are ACD. 
	For this we can replace the elliptic surface $E(n)$ without multiple 
	fibers by one with multiple fibers obtained by logarithmic transformation 
	as in the proof of Proposition~\ref{p:wedge}. As all elliptic surfaces 
	are ACD, and the logarithmic transformations can be assumed to have been 
	made inside $E(n-1)$ in a splitting $E(n)=E(n-1)\cup_{T^{2}} E(1)$, 
	all the resulting manifolds will be ACD by the same argument as above.
\end{proof}

\section{The main theorems}\label{s:proof}

The following theorem is very close to various results proved by 
Mandelbaum and Moishezon, and will be proved using their technique,
but it does not appear explicitly in their papers~\cite{Ma,Msurv,MM}. 
The case of complete intersections does appear there, and it is pointed out 
that they are branched covers, but the latter are not treated in 
complete generality.

\begin{thm}\label{salvetti-acd}
Iterated branched covers of the projective plane are almost completely 
decomposable.
\end{thm}
\begin{proof}
    The proof is by induction on the number of iterations. To begin, 
    note that the cyclic cover of degree $d$ of the complex projective plane 
    branched in a smooth curve of degree $d\cdot m$ is ACD by applying 
    Theorem 2.9 of~\cite{Msurv} to the Veronese embedding of $\C P^{2}$ 
    given by the monomials of degree $m$. 
    
    It remains to show that if $f\colon X\map\C P^2$ is an iterated branched 
    cover and $Y$ is the cyclic branched cover of $X$ branched over 
    $f^{-1}(C_d)$, where $C_d$ is a general plane curve of degree $d\cdot m$, 
    then $Y$ is ACD if $X$ is.
    
    Keeping in mind that $f^{-1}(C_{d})$ is ample, one can find closely 
    related statements in~\cite{Msurv}, yet we cannot rely on them directly. 
    In Theorem 2.9 of~\cite{Msurv} the branch locus is assumed to be very 
    ample, whereas Theorem~2.14 refers to the homology class of the branching 
    locus and is not applicable to a given representative. Nevertheless we 
    follow the line of argument of Mandelbaum and Moishezon, see in particular 
    Theorem~4.2 in~\cite{Ma} or Theorems~4.1 and 4.2 in~\cite{MM}.
    
    For each $k\geq 0$ choose a general section $s_k$ of the degree $km$
    line bundle on the plane and let $C_k$ be its vanishing locus. Then 
    $f^{-1}(C_k)$ is a smooth curve in $X$ given by the equation 
    $f^\ast s_k=0$. Consider the line bundle $L=f^\ast\mathcal O(m)$ and its 
    compactification $W=\PP(L\oplus\mathcal O_X)$. If $p\colon W\map X$ is the 
    bundle projection and $W_\infty=W\weg L$ is the section at infinity then 
    the line bundle $E=p^\ast(L)(W_\infty)$ admits a tautological section $y$ 
    without zeroes at infinity.  A $k$-sheeted cyclic covering $Y_k$ of $X$ 
    branched over $f^{-1}(C_k)$ is described in $W$ as the vanishing locus of 
    the section $t_k=y^k-(fp)^\ast s_k$ of $E^k$. Consider the pencil in 
    $|E^d|$ generated by $t_d$ and $t_1\cdot t_{d-1}$. The general member 
    of the pencil is smooth and diffeomorphic to $Y_d$. The special member 
    given by the vanishing of $t_1t_{d-1}$ is the union of two smooth surfaces 
    $Y_1\cong X$ and $Y_{d-1}$ intersecting in a curve which is given on $Y_1$ 
    by the equation $y^{d-1}=(fp)^\ast s_{d-1}$. The obvious isotopy 
    from $Y_1$ to $X$ (embedded into $W$ as the zero section) transforms this 
    curve into the curve on $X$ given by $0=f^\ast s_{d-1}$, 
    i.~e.~$f^{-1}(C_{d-1})$.  If this is a sphere then its preimage in 
    $Y_{d-1}$ is a sphere with positive self-intersection.
    In this case $Y_{d-1}$ is rational hence completely decomposable.

On the other hand, if the curve $f^{-1}(C_{d-1})$ has genus at least $1$ we
consider the fibration of the blowup of $W$ along $(Y_1\cup Y_{d-1})\cap Y_d$
over $\C P^{1}$ given by our pencil.  By~\cite{MM} Corollary 2.7, $Y_d$ is the
sum of a blowup of $Y_1\cong X$ and of $Y_{d-1}$ along $f^{-1}(C_{d-1})$.
By~\cite{MM} Theorem 2.8 (2) it follows that $Y_d\#\C P^2$ is diffeomorphic to
the connected sum of $X$ and $Y_{d-1}$ together with $k>0$ copies of $\C P^2$
and $\overline{\C P^2}$.  Since $X$ is almost completely decomposable
the result follows by induction.
\end{proof}

\begin{rem}
    The same argument applies to the iterated branched 
    covers of a quadric considered by Moishezon~\cite{Mjdg}.
\end{rem}

Next we exhibit manifolds covering a large geographical area 
satisfying the Hitchin-Thorpe inequality~\eqref{eq:HT}, but which have 
infinitely many smooth structures not supporting Einstein metrics.
\begin{thm}\label{t:non-Einstein}
    For every $\epsilon >0$, there is a constant $c(\epsilon)>0$ such 
    that every lattice point $(x,y)$ with $y\geq 0$ satisfying
    $$
    y\leq (6-\epsilon)x-c(\epsilon)
    $$
    is realized by the Chern invariants $(\chi,c_{1}^{2})$ of infinitely 
    many pairwise non-diffeomorphic simply connected almost completely 
    decomposable symplectic manifolds which do not admit Einstein 
    metrics.
    \end{thm}
\begin{proof}
    We consider the manifolds $Y(l,k,r,n)$ in the proof of Theorem~\ref{t:geo} 
    above. They are all symplectic, and so have non-trivial Seiberg-Witten 
    invariants. Therefore~\cite{lebrun,K}, if such a manifold is blown up 
    sufficiently often, the blowup cannot admit any Einstein metric. 
    According to Theorem~3.3 of LeBrun~\cite{lebrun3}, 
    $\frac{1}{3}c_{1}^{2}(Y(l,k,r,n))$ many blowups suffice. Thus, 
    these manifolds cover the claimed area. The infinitely many 
    distinct smooth structures on each remain distinct under blowing 
    up, see for example~\cite{FS,KMT}. Clearly the ACD property is 
    preserved by the blowups.
    \end{proof}

    We can now combine the results proved so far in order to prove the 
    following more detailed version of Theorem~\ref{t:main}.
\begin{thm}\label{t:main2}
    For every natural number $k$ there are simply connected topological 
    $4$-manifolds $M_{k}$ which have at least $k$ distinct smooth structures 
    $Z_{k}^{i}$ supporting Einstein metrics, and also have infinitely many 
    distinct smooth structures $X_{k}^{j}$ not supporting Einstein metrics. 
    
    The $Z_{k}^{i}$ and $X_{k}^{j}$ can be chosen symplectic and almost 
    completely decomposable. For every fixed $k$, the ratios $c_{1}^{2}/\chi$ 
    of the Chern invariants of such examples are dense in the interval $[4,6]$.
    \end{thm}
\begin{proof}
    We consider certain simply connected symplectic manifolds 
    which are non-spin and have the same Chern invariants. Thus they are 
    homeomorphic by Freedman's classification~\cite{freed}.
    
    The $Z_{k}^{i}$ are the iterated branched covers of the projective 
    plane constructed in Theorem~\ref{ratio-4.8}. By Corollary~\ref{c:geography}, 
    the ratios of the Chern invariants of such examples are dense in the interval 
    $[4,8]$. By Theorem~\ref{salvetti-acd}, these manifolds are ACD. 
    As they are K\"ahler with ample canonical bundle, the solution of 
    the Calabi conjecture due to Aubin~\cite{A} and Yau~\cite{Y} shows 
    that they carry Einstein metrics.
    
    Bringing down the upper bound for the slope to $6$ allows us to use 
    manifolds from Theorem~\ref{t:non-Einstein} having appropriate 
    Chern invariants for the $X_{k}^{j}$. These are ACD by construction 
    and do not carry Einstein metrics. 
    
    We already noted in Corollary~\ref{c:geography} that the $Z_{k}^{i}$ 
    are pairwise non-diffeomorphic by Seiberg--Witten theory~\cite{Wi}. 
    The $X_{k}^{j}$ are obtained by blowing up distinct smooth structures 
    distinguished by their Seiberg-Witten invariants, and so they are 
    also distinct because of the blowup formula~\cite{FS,KMT}. Clearly 
    no $Z_{k}^{i}$ can be diffeomorphic to a $X_{k}^{j}$, as the former 
    admit Einstein metrics and the latter do not. (Also, the former are 
    irreducible~\cite{Bourbaki}, and the latter are not.)
\end{proof}
\begin{rem}\label{r:FS}
The manifolds $M_{k}$ have another infinite sequence of smooth structures, 
which are very likely distinct from the $Z_{k}^{i}$ and the $X_{k}^{j}$. 
Fintushel--Stern~\cite{FScusp} showed that one can perform cusp surgery on 
a torus in any iterated branched cover of the plane to construct infinitely 
many distinct smooth structures with non-trivial Donaldson invariants. It 
seems that these are irreducible, and therefore distinct from the $X_{k}^{j}$. 
On the other hand they are not complex, and therefore distinct from the 
$Z_{k}^{i}$. Whether they are ACD or admit Einstein metrics is not known.
\end{rem}

Theorems~\ref{t:geo} and~\ref{t:non-Einstein} also lead to the following 
more general existence result for smooth structures not supporting Einstein 
metrics.
\begin{thm}\label{t:pq}
    For every $\epsilon>0$ there is a constant $c(\epsilon)>0$ such 
    that the connected sum $p\C P^{2}\# q \overline{\C P^{2}}$ has 
    infinitely many smooth structures not admitting Einstein metrics 
    for every large enough $p\not\equiv 0\pmod 8$ and 
    $q\geq (2+\epsilon)p+c(\epsilon)$.
    \end{thm}
    \begin{proof}
	For odd $p$, this was already proved in Theorem~\ref{t:non-Einstein}.
	
	For even $p$, we are in a situation where the numerical 
	Seiberg--Witten invariants must vanish. Therefore, to obtain an 
	obstruction to the existence of Einstein metrics one considers the 
	refined Seiberg--Witten invariants of Bauer--Furuta~\cite{BF} in the 
	context of stable homotopy theory. Using this approach, 
	Ishida--LeBrun~\cite{IL} showed that a connected sum 
	$X_{1}\# X_{2}\# k\overline{\C P^{2}}$, where the $X_{i}$ are simply 
	connected symplectic four-manifolds with $b_{2}^{+}\equiv 3 \pmod 4$, 
	does not admit Einstein metrics if 
	$k\geq\frac{1}{3}(c_{1}^{2}(X_{1})+c_{1}^{2}(X_{2}))-4$. Applying 
	this to the case where $X_{1}$ are the manifolds from 
	Theorem~\ref{t:geo} with $b_{2}^{+}\equiv 3 \pmod 4$ and $X_{2}$ 
	is the $K3$ surface, proves the claim of the Theorem for 
	$p\equiv 2\pmod 4$. We just have to see that the connected sum with 
	$K3$ does not collapse the infinitely many smooth structures on 
	$X_{1}\# k\overline{\C P^{2}}$. These smooth structures were 
        constructed by 
	logarithmic transformation on an elliptic building block in $X_{1}$. 
	As we increase the multiplicity of the logarithmic transformation, 
	we find that there are more and more Seiberg--Witten basic classes 
	whose numerical Seiberg--Witten invariants are $\pm 1$, see 
	Fintushel--Stern~\cite{FSrat}, Theorem~8.7. By the result of
        Bauer~\cite{B}, these basic classes give rise to monopole classes 
        in the sense of~\cite{Kr} on $X_{1}\# X_{2}\# k\overline{\C P^{2}}$. 
        As the expected dimension of the Seiberg--Witten moduli space is 
        positive for all these monopole classes, each smooth structure has
        at most finitely many such classes. This shows that we have an 
        infinite set of smooth structures\footnote{See~\cite{K3} for 
        more details and elaborations on this argument.}.

	
	It remains to deal with the case $p\equiv 0\pmod 4$. The above argument 
	generalizes to the case of connected sums of $4$ symplectic manifolds 
	$X_{i}$ with $b_{2}^{+}\equiv 3 \pmod 4$ as long as the resulting 
	manifold $X_{1}\#\ldots\# X_{4}\# k\overline{\C P^{2}}$ has $b_{2}^{+}$ 
	not divisible by $8$ and 
	$k\geq\frac{1}{3}(c_{1}^{2}(X_{1})+\ldots+c_{1}^{2}(X_{4}))-12$. 
	This was noted by Ishida--LeBrun in~\cite{IL2}, using~\cite{B}. We apply 
	it here taking for $X_{1}$ the manifolds from Theorem~\ref{t:geo} with 
	$b_{2}^{+}\equiv 3 \pmod 8$, and taking the $K3$ surface for $X_{2}$, 
	$X_{3}$ and $X_{4}$. This proves the claim of the Theorem for 
	$p\equiv 4\pmod 8$.
	\end{proof}
\begin{rem}
Theorem~\ref{t:pq} should be compared to Theorems 11 and 12 of 
Ishida--LeBrun~\cite{IL}, which give much weaker statements in the same 
direction. Namely, if $p$ is odd they assumed $p\equiv 1 \pmod 4$ and 
$q>\frac{7}{3}p+12$, which is more restrictive than 
$q\geq (2+\epsilon)p+c(\epsilon)$ for almost all $p$ whenever 
$\epsilon<\frac{1}{3}$. The unknown constant $c(\epsilon)$ only appears 
in our statement because we constructed all manifolds to be ACD. If we 
gave up this constraint, we could make the constant explicit. However, 
our method of proof, and the smooth structures under consideration, are 
very different. In our proof, for odd $p$ the smooth structures in question 
support symplectic forms, and, therefore~\cite{Bourbaki}, can not 
decompose as smooth connected sums except for the splitting off of copies 
of $\overline{\C P^{2}}$. The non-existence of Einstein metrics is detected 
by the numerical Seiberg--Witten invariants. The smooth structures 
discussed by Ishida--LeBrun~\cite{IL} are smooth connected sums where each 
summand has positive $b_{2}^{+}$, and so in particular they cannot support 
symplectic forms. As the numerical Seiberg--Witten invariants of these 
smooth structures vanish, the non-existence of Einstein metrics can only be 
detected using the stable homotopy refinement~\cite{BF} of the Seiberg--Witten 
invariants. 

For even $p$, Ishida--LeBrun~\cite{IL} assumed $p\equiv 2 \pmod 4$ and 
$q>\frac{7}{3}p+16$. In this case our proof is similar to 
theirs---in fact we use the main result of their paper. Our improvement 
is due to the fact that our Theorem~\ref{t:geo} above gives us more symplectic 
manifolds we can use as connected summands, whereas Ishida and LeBrun used 
only certain manifolds constructed by Gompf~\cite{Go1} with smaller slope of 
their Chern invariants.
\end{rem}

\section{Further examples}\label{s:ex}

Since the proof of Theorem~\ref{t:GT} in~\cite{K}, other examples of 
manifolds with a smooth structure supporting an Einstein metric and 
one or several without an Einstein metric have appeared, and some 
attempts have been made to give examples with smallish homology, 
compare~\cite{lebrun3,IL}. Here are the ultimate examples, whose 
second Betti number is a fraction of that of the smallest previously 
known examples.

\begin{prop}\label{t:small}
    The manifolds $3\C P^{2}\# 17\overline{\C P^{2}}$ and 
    $3\C P^{2}\# 18\overline{\C P^{2}}$ each have a smooth structure 
    supporting an Einstein metric, and infinitely many smooth 
    structures not supporting Einstein metrics.
    \end{prop}
\begin{proof}
    Take a double cover of $\C P^{2}$ branched in the union of two smooth 
    cubics in general position. This gives a singular $K3$ surface with 
    $9$ nodes. Now take a further double cover of this singular $K3$ 
    surface branched in the nodes and in the preimage of a line. This gives 
    a simply connected smooth algebraic surface $S$ with ample canonical 
    bundle, whose numerical invariants are $c_{1}^{2}(S)=1$, $\chi(S)=2$, 
    compare~\cite{Ca0}. It is homeomorphic to $3\C P^{2}\# 18\overline{\C P^{2}}$ 
    by Freedman's classification~\cite{freed}. By the results of Aubin~\cite{A} 
    and Yau~\cite{Y} it admits a K\"ahler--Einstein metric.
    
    One way to obtain a smooth manifold homeomorphic to $S$ which cannot 
    admit an Einstein metric is to take a simply connected algebraic surface 
    $S'$ with $c_{1}^{2}(S')=2$, $\chi(S')=2$, and blow it up once. Such $S'$ 
    exist, compare Catanese--Debarre~\cite{CD} and the discussion below.
    According to LeBrun~\cite{lebrun3}, the blowup of $S'$ does not admit 
    any Einstein metric. Another possibility is to take a symplectic manifold 
    homeomorphic to $3\C P^{2}\# n\overline{\C P^{2}}$ for some $n\leq 16$, 
    and blow it up until it becomes homeomorphic to $S$. Such manifolds have 
    been constructed by Gompf~\cite{Go1} and D.~B.~Park~\cite{DP}. There are in 
    fact infinite sets of smooth structures on them supporting symplectic forms, 
    compare~\cite{P3}. These remain distinct under blowing up, and the blowups 
    have no Einstein metrics by the result of~\cite{lebrun,K}, say. This proves 
    the claim for $3\C P^{2}\# 18\overline{\C P^{2}}$.
    
    To obtain an algebraic surface homeomorphic to 
    $3\C P^{2}\# 17\overline{\C P^{2}}$ which has ample canonical 
    bundle one can proceed as follows. Take a double cover of $\C 
    P^{1}\times \C P^{1}$ branched in the union of two smooth curves 
    of bidegrees $(3,1)$ and $(1,3)$ respectively. Then take a further 
    double cover branched in the nodes of the first covering and the 
    preimage of a smooth curve of bidegree $(1,1)$ in general position 
    with respect to the other two curves. The resulting smooth surface $S'$ 
    has all the desired properties, compare~\cite{CD}. If we start with one 
    of the symplectic manifolds with $c_{1}^{2}\geq 6$ constructed by 
    D.~B.~Park~\cite{DP}, then there are infinitely many smooth 
    structures on it which remain distinct under blowing up, and the 
    blowups homeomorphic to $S'$ do not admit any Einstein metrics 
    by~\cite{lebrun,K}.
    \end{proof}

    Catanese~\cite{Ca0} proved that all the algebraic surfaces $S$ 
homeomorphic to $3\C P^{2}\# 18\overline{\C P^{2}}$ are diffeomorphic 
to each other. In~\cite{CD} it is conjectured that the same is true 
for surfaces homeomorphic to $3\C P^{2}\# 17\overline{\C P^{2}}$. 
Thus, one has to take larger examples to obtain multiple smooth 
structures with Einstein metrics\footnote{It is known that 
$\C P^{2}\# 8\overline{\C P^{2}}$ has two distinct smooth structures 
supporting Einstein metrics, see~\cite{CL,K}, but it is not known 
whether this manifold has a smooth structure without an Einstein 
metric.}. While our proof of Theorem~\ref{t:main} can be made effective, 
in practice the manifolds one obtains will have huge homology. Nevertheless, 
concrete examples can be given quite easily.
\begin{ex}
    For any integer $k\geq 0$ let $Z_{2}^{1}$ be a smooth hypersurface of 
    bidegree $(5+k,6)$ in $\C P^{1}\times\C P^{2}$, and let $Z_{2}^{2}$ be a 
    smooth complete intersection of two hypersurfaces of bidegrees $(2,1)$ 
    and $(1+k,6)$ in $\C P^{1}\times\C P^{3}$. Both have 
    $$
    c_{1}^{2}=9(17+5k) \ ,
    $$ 
    $$
    \chi=41+10k \ . 
    $$
    The divisibility of the canonical class is $gcd\{k+3,3\}=gcd\{k,3\}$ for 
    $Z_{2}^{1}$ and $gcd\{k+1,3\}$ for $Z_{2}^{2}$. Thus they are both non-spin 
    and are homeomorphic for each $k$. If $k\equiv 1\pmod 3$, then they both 
    have divisibility $=1$, otherwise they have different divisibilities, and 
    can therefore not be diffeomorphic. 
  
These surfaces have ample canonical bundles, and therefore~\cite{A,Y} support 
K\"ahler--Einstein metrics. Their ratio $c_{1}^{2}/\chi$ is $<4.5$, and thus they 
are well within the range where Theorem~\ref{t:non-Einstein} shows that there are 
infinitely many homeomorphic smooth manifolds without Einstein metrics.
       \end{ex}
       
Finally, as the manifolds discussed so far are all non-spin, it is 
worthwhile to point out the following:
\begin{thm}\label{t:evenSalvetti}
    For every $k$ there is a topological spin four--manifold $M_{k}$ 
    admitting at least $k$ distinct smooth structures which support Einstein 
    metrics, and which are almost completely decomposable. The ratios 
    $\vert\sigma\vert/e$ of such manifolds are dense in the interval 
    $[\tfrac{1}{3},\tfrac{1}{2}]$.
    \end{thm}
\begin{proof}
    This follows from the spin case of Corollary~\ref{c:geography} 
    together with the existence result for K\"ahler--Einstein metrics 
    due to Aubin~\cite{A} and Yau~\cite{Y}. The ACD property was proved in 
    Theorem~\ref{salvetti-acd}.
    \end{proof}
\begin{rem}\label{r:cat}
    Such examples are also provided by Catanese's $k$-tuples of homeomorphic 
    spin surfaces with different divisibilities of their canonical classes 
    constructed as bidouble covers of a quadric~\cite{Ca2}. For almost all 
    choices of the parameters, those surfaces have ample canonical bundles. 
    However, the spread of their numerical invariants is probably more 
    restricted than in our examples based on Salvetti's 
    construction~\cite{Salve}.
\end{rem}
  
\bibliographystyle{amsplain}

\bigskip

\end{document}